# Grade of the law of distribution of the vector of the digital twin of the enterprise

*Sergei* Masaev[1*], *Ivan* Oreshnikov[1], *Nikita* Ivanitskiy[1], and *Valentina* Vingert[2]

[1]Reshetnev Siberian State University of Science and Technology, 31, Krasnoyarsky Rabochy Avenue, Krasnoyarsk, 660037, Russian Federation
[2]Siberian Federal University, 79, prospect Svobodnyj, Krasnoyarsk, 660041, Russian Federation

**Abstract.** Digital transformation acts as a main factor in the development and competitiveness of Japanese companies. In this context, the digital copy management algorithm based on the P2M (Project to Method) method plays a significant role in increasing efficiency and reducing costs through the implementation of digital technologies. This research aims to develop and implement an algorithm for managing a digital copy of a project using the P2M method using the Kolmogorov test to evaluate the effectiveness of various processes and determine the best ways to improve productivity.

## 1 Introduction

Effective control is becoming a key success factor in enterprises with an ever-increasing complexity of the external environment in the modern world. The Project Management Association of Japan (PMAJ) P2M method is a comprehensive approach to managing innovation projects. This relies on the interaction of three main components: complexity, value and resistance [1].

There is no formalized approach to estimating its effectiveness in managing an enterprise as a dynamic system in Japan. Existing modeling methods do not take into account the specifics of P2M, which makes it difficult to identify and interpret patterns associated with its influence [1].

The purpose of this research is to estimate the state of an economic object, considered as a multidimensional dynamic system when applying the P2M management method under conditions of uncertainty of external parameters. The Kolmagorov test used will allow us to analyze and evaluate the effectiveness of implementing this control method to achieve optimal results at the enterprise [2-5].

## 2 Method

The beginning of the activity of a timber industry enterprise was chosen as a digital project (object of research). The project is located in the North Yenisei region. The company harvests 800 thousand cubic meters of roundwood, which is then transported along the

---

* Corresponding author: faberi@list.ru

Yenisei River on barges from June to September. Harvested raw materials are used to produce advanced wood processing products, including floorboards, laminated veneer lumber, eurolining and other products [6-9].

The company intends to increase production volumes using bank loans and tax incentives over the next 1.5 years. The success of this strategy depends on market development and the strategy chosen by the company [10]. The enterprise is implementing a management system based on the P2M method in this dynamic situation [1]. The enterprise implements a control loop based on the P2M method in the context of this dynamic situation [1]. This approach involves integrated management of various aspects of the project, including:

- Project strategy management: defining long-term goals, strategic objectives and ways to achieve them;
- Project financing management: planning, control and optimization of financial resources [2];
- Project systems management: development, implementation and maintenance of the necessary systems and technologies [3];
- Organizational project management: formation of a highly productive team and improvement of interaction between its members;
- Managing project goals: clearly defining the goal, tracking the progress of implementation and achieving the set goal of the project;
- Project resource management: planning, distribution of resources (material, financial, human) [4];
- Application of modern project information technologies for effective management and exchange of data within the project [8];
- Formation of relationships in the project, ensuring achievement of the project's objectives;
- Management, monitoring and reduction of project implementation costs [6];
- Ensuring open data exchange between project participants [7].

## 3 P2M method (Japan)

The P2M (Project and Program Management) method appeared in Japan in November 2001. The main essence on which this method is based is the control of projects and software products. The principles of the P2M method are analysis of the project from the point of view of creating new value that brings benefits to the client [8, 9].

Basic principles:
- Planning: setting the main goals of the project, distributing tasks, determining deadlines, and assigning responsible persons;
- Development: the project is divided into smaller tasks (subtasks), which are distributed among team members;
- Control: tracking progress at each stage of the project and identifying problems;
- Management: managing risks, resources and project communications.

## 4 Kolmogorov-Smirnov test

The Kolmogorov-Smirnov test is a parameter-independent statistical test used to test the difference between two probability distributions. It was developed by Andrei Kolmogorov and Nikolai Smirnov in 1933. Let's say we want to find out whether there are differences between the distributions of random variables *X* and *Y*. After we compare how much the actual distributions of these variables differ from their expected (theoretical) distributions [5].

The test statistic $D_n$ is defined as the maximum difference between these functions and is expressed by the formula:

$$D_n = sup_x |Fn(x) - F(x)| \qquad (1)$$

The result of testing the hypothesis is two statements:
- If the p-value is less than a certain significance level, then the hypothesis of equality of distributions can be rejected ($\boldsymbol{H_0}$);
- If the *p*-value is greater than the significance level, then the distributions are different ($\boldsymbol{H_1}$) [5].

## 5 Vector testing process

Based on the performed Kolmogorov-Smirnov test to estimate the normality of data distribution, the following conclusions can be drawn [5].

The null hypothesis (*H₀*) of normal distribution was tested using the Kolmogorov-Smirnov test. The formula for the Kolmogorov-Smirnov statistics is (Fig. 1, Figures 1 and 2):

$$D_n = sup_x |Fn(x) - F(x)| \qquad (1)$$

— empirical distribution function, *F(x)* — cumulative normal distribution function. The test statistic is compared with the critical value for a given significance level *α* [5].

**Table 1.** Initial values (data presented partially).

| Variables | Initial data | |
|---|---|---|
| | without P2M | by P2M |
| x1 | 34654.53 | 34654.53 |
| | 34919.47 | 34919.47 |
| | 88541.03 | 86456.58 |
| | 88661.97 | 86577.52 |
| | 60466.67 | 58382.23 |
| | 162752.50 | 160668.06 |
| | 159851.27 | 157816.22 |
| | 156642.61 | 154607.56 |
| | 160696.10 | 158661.04 |
| | 144506.44 | 142471.38 |
| x2 | 51017.81 | 51017.81 |
| | 51765.31 | 51765.31 |
| | 92337.28 | 92337.28 |
| | 91291.56 | 91291.56 |
| | 139128.51 | 139128.51 |
| | 48591.03 | 48590.50 |
| | 42175.84 | 42291.09 |
| | 38758.30 | 38873.55 |
| | 47978.81 | 48093.54 |
| | 39456.58 | 39571.83 |
| x3 | 3559.08 | 3032.24 |
| | 6684.59 | 6299.49 |
| | 728767.56 | 729474.57 |
| | 9672.61 | 9858.69 |
| | 9225.86 | 9411.93 |
| | 9714.77 | 9900.85 |
| | 10292.01 | 10675.60 |
| | 9311.71 | 9668.10 |
| | 9919.97 | 10306.29 |
| | 9466.93 | 9823.18 |
| x4 | 89231.42 | 88704.58 |
| | 93369.37 | 92984.27 |
| | 909645.88 | 908268.44 |
| | 189626.14 | 187727.77 |

|  | 208821.04 | 206922.67 |
|  | 221058.30 | 219159.40 |
|  | 212319.12 | 210782.90 |
|  | 204712.62 | 203149.20 |
|  | 218594.89 | 217060.87 |
|  | 193429.95 | 191866.39 |

## 6 Test result

**Table 2.** Test results (partially presented).

| Variables | D - statistics | P-value | Conclusion |
|---|---|---|---|
| X1 | 0.190097 | 0.022588 | The null hypothesis is rejected |
| X2 | 0.335311 | 0.000002 | The null hypothesis is rejected |
| X3 | 0.189035 | 0.023738 | The null hypothesis is rejected |
| X4 | 0.339469 | 0.000001 | The null hypothesis is rejected |

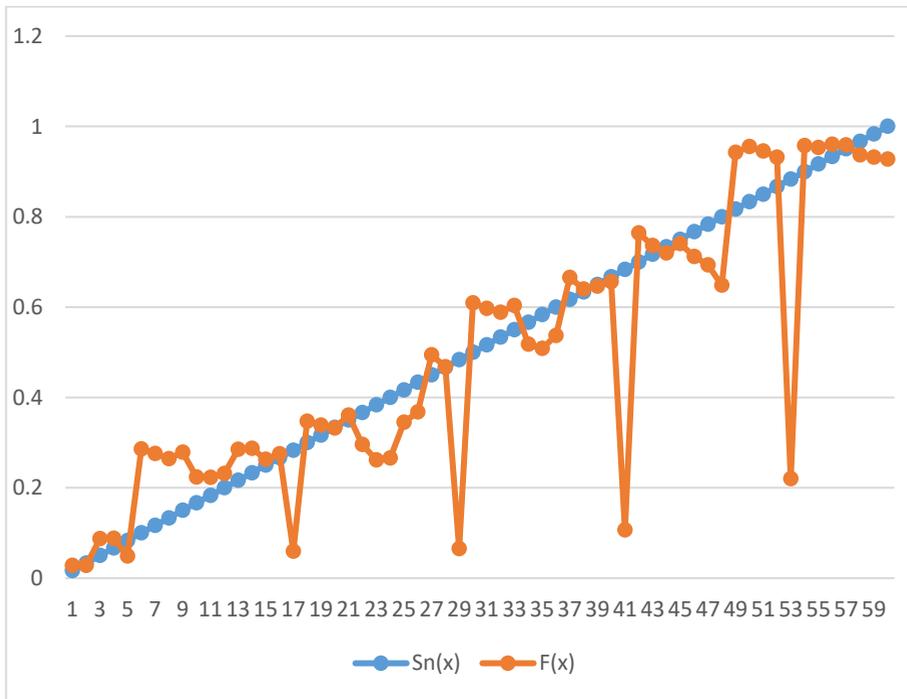

**Fig. 1.** Evaluation result of the Kolmogorov-Smirnov test.

## 7 Discussion of results

The calculations and analysis revealed significant deviations from the normal distribution for data *X1, X2, X3* and *X4* by the Kolmogorov-Smirnov test. The *D* statistic was 0.190097 with a *P*-value of 0.022588, indicating a noticeable deviation from the normal distribution for variable *X1*. Variable *X2* shows the largest deviation with $D = 0.335311$ and *P*-value = 0.000002. Statistics for variable *X3* showed a value of $D = 0.189035$ with a *P*-value of 0.023738 so confirms significant deviations from the normal distribution. Variable *X4* has a statistic of $D = 0.339469$ and *P*-value = 0.000001 indicating a significant deviation from the normal distribution.

The identified deviations from the normal distribution require a revision of current methods and approaches to data analysis. Understanding these variations will help develop more accurate and reliable models for estimating and predicting manufacturing processes. These adjustments will help improve the accuracy and reliability of the data obtained as well as improves the overall efficiency of processes.

## 8 Conclusion

The results of the Kolmogorov-Smirnov test to estimate the normality of data distribution showed that the hypothesis of normal distribution is rejected for the following data: *X1, X2, X3* and *X4*. This means that the distribution of data in these columns is significantly different from the normal distribution at the 0.05 significance level. For the remaining data, the hypothesis of a normal distribution is also rejected, indicating that the data deviate from the normal distribution in these columns.

The findings provide important information for further data analysis and selection of adequate analysis methods based on the Kolmogorov-Smirnov test. The purpose of the research has been achieved.